\documentclass[a4paper, oneside, notitlepage, 10pt]{amsart}

\usepackage{amsmath}
\usepackage{amsfonts}
\usepackage{amssymb}
\usepackage{amsthm}
\usepackage{graphicx}
\usepackage{mathrsfs}

\RequirePackage{hyperref}

\usepackage{verbatim}
\usepackage[T1]{fontenc}        
\usepackage[english]{babel}                              
\usepackage[latin1]{inputenc}  

\usepackage{amsmath,amsfonts,amssymb,amsthm,graphicx,mathrsfs}

\usepackage{verbatim}        
\usepackage[english]{babel}                              
\usepackage[latin1]{inputenc}  
\usepackage{fontenc}

\DeclareMathOperator{\Xb}{\overline{X}(\infty)}
\DeclareMathOperator{\Xg}{\partial \overline{X}}
\DeclareMathOperator{\X}{\overline{X}}
\DeclareMathOperator{\Isom}{Isom}
\DeclareMathOperator{\CAT}{CAT}

\newcommand{\sra}{\rightarrow}                  % ->

\newcommand{\N}{\mathbb{N}_0}                 % Natural numbers
\newcommand{\Z}{\mathbb{Z}}                   % Integers
\newcommand{\R}{\mathbb{R}}                  % Real numbers
\newcommand{\mas}[1]{\{ #1 \}}                   % { }
\newcommand{\pan}[1]{(#1)}                        % ( )
\newcommand{\hak}[1]{[#1]}                         % [ ]

\newcommand{\phib}{\varphi}
\newcommand{\eps}{\varepsilon}

\newcommand{\Pb}{\mathbb{P}}
\newcommand{\Pbt}{\tilde{\mathbb{P}}}

\renewcommand{\H}{\mathcal{H}}

\theoremstyle{plain} 
\newtheorem{thm}{Theorem}
\newtheorem{lem}{Lemma}[section]

\theoremstyle{definition}

\newtheorem{rema}{Remark}[section]
\newtheorem{exam}{Example}[section]
\newtheorem*{ass}{Basic Assumption (BA)}

\theoremstyle{plain}

\begin{document}

\title{Central Limit Theorems for Gromov Hyperbolic Groups}
\author{Michael Bj\"orklund}
\email{mickebj@math.kth.se}
%\date{}
\maketitle

\begin{abstract} 
In this paper we study asymptotic properties of symmetric and non-degenerate random walks on transient hyperbolic groups. We prove a central limit theorem 
and a law of iterated logarithm for the drift of a random walk, extending previous results by S. Sawyer and T. Steger and F. Ledrappier for certain $ \CAT(-1)$--groups. 
The proofs use a result by A. Ancona on the identification of the Martin boundary of a hyperbolic group with its Gromov boundary. We also give a new interpretation, in 
terms of Hilbert metrics, of the Green metric, first introduced by S. Brofferio and S. Blach\`ere. 
\end{abstract}

\section{Introduction}

Let $ X_1 , X_2 , \ldots $ be a sequence of independent and identically distributed random variables on the real line, such that $ \mathbb{E}\hak{X_k} = 0 $ and $ 0 < \mathbb{E}\hak{X_k^2} < \infty $. 
The central limit theorem tells us that the expression
\[
\frac{X_1 + \ldots X_n}{\sqrt{n}} 
\]
converges in distribution to a non-degenerate Gaussian variable. Bellman \cite{Be54} and Furstenberg and Kesten \cite{FuKe60} initiated the study of non-commutative analogues of the 
central limit theorem. In \cite{FuKe60}, Furstenberg and Kesten proved that if $ X_1 , X_2 \ldots $ is a sequence of bounded, independent and identically distributed random matrices of a 
fixed size $ d $, then, if $ Y_n = X_1 \cdots X_n $, and under some technical assumptions, the sequence
\[
\frac{\log{(Y_n)}_{ij}-\mathbb{E}\hak{\log{(Y_n)}_{ij}}}{\sqrt{n}}
\] 
converges weakly to gaussian variables for all $ i , j = 1 , \ldots , d $.

In this paper, we consider the following related problem: 

Let $ ( X , d ) $ be a metric space, and $ g_1 , g_2 , \ldots $ is a sequnce of independent and identically distributed random variables taking 
values in the isomety group of $ ( X , d ) $. Pick a basepoint $ x_0 \in X $ and let $ A_n = d(g_1 \cdots g_n x_0,x_0) $.  We assume that $A_1$ has finite first moment. 
Then, by Kingman's subadditive ergodic theorem \cite{Ki73}, the limit
\[
A = \lim_{n \sra \infty} \frac{A_n}{n}
\]
exists almost surely, and is constant almost everywhere. We may ask about the rate of convergence to $ A $. This question heavily depends on the metric space $ ( X , d ) $, and no general results are known. In \cite{SaSt87},
Sawyer and Steger studied the case when $ ( X , d ) $ is the Cayley graph of a free group $ F_d $ on $ d $ standard generators: Let $ g_1 , g_2, \ldots $ be i.i.d. random variables taking values in $F_d$, and distributed
according to a probability measure $ \mu $ on $ F_d $, which is required to satisfy some technical moment conditions, then $ A > 0 $, and
\[
\frac{d(g_1 \cdots g_n,e)-nA}{\sqrt{n}}
\] 
converges to a non-degenerate Gaussian distribution on $ \R $. Their proof relies heavily on the properties of certain analytic extensions of functions of the Green kernel of the measure $ \mu $. A geometric proof 
was later given by F. Ledrappier in \cite{Le01} and a nice generalization to random walks on trees with finitely many cone types is proved in \cite{NaWo02}. The present paper is concerned with a partial generalization of 
these results to random walks on a general Gromov hyperbolic group $ \Gamma $. 

Let $ p $ be the uniform measure on a finite set $ S $ of generators of $ \Gamma $. If we exclude the case of finite groups and finite extensions of $ \Z $ and $ \Z^2$, the Green kernel
\[
G(x,y) = \sum_{n=0}^{\infty} p^{*n}(y^{-1}x) 
\]
is everywhere positive and finite. It turns out that $ d(x,y) = \log{\frac{G(e,e)}{G(x,y)}} $ defines a $ \Gamma $-invariant metric on $ \Gamma $, and roughly equivalent to the word-metric with respect to 
$ S $, in the sense that there are positive constants $C$ and $b$ such that
\[
\frac{1}{C} d_S(x,y) - b \leq d(x,y) \leq Cd_S(x,y) + b, \quad \forall \: x,y \in \Gamma. 
\]
We will refer to $ d $ as the Green metric on the Cayley graph of $ (\Gamma,S) $. This metric was first introduced in \cite{BlBr07}, but has been studied quite a lot in disguise, in connection with Martin 
boundaries of groups \cite{An87}, \cite{An90}. Using this metric, we can formulate the following theorem

\begin{thm} \label{thm:main}
Let $ \Gamma $ be a Gromov hyperbolic group, which is not finite or virtually $ \Z $. If the support of $ \mu $ generates $ \Gamma $ and if there is a $ \beta > 0 $, such that $ \int_{\Gamma} e^{\beta d(g,e)} \, d\mu(g) < \infty $, then if
$ g_1 , g_2 , \ldots $ is a sequence of independent $ \mu$--distributed random variables, then $ A > 0 $ and, 
\[
\frac{d(g_1 \cdots g_n , e)-nA}{\sqrt{n}} 
\]
converges to a non-degenerate centered Gaussian variable, with variance $ \sigma > 0 $, and furthermore, 
\[
\limsup_{n \sra \infty} \frac{d(g_1 \cdots g_n,e)-nA}{\sqrt{n\log{\log{n}}}} = \sigma
\]
almost surely. 
\end{thm}
Note that this theorem does not imply the central limit theorem for the drift with respect any word--metric on $\Gamma$. 

The reason for introducing the Green metric is geometric. With respect to the Green metric, the horofunction boundary of $ ( \Gamma , d ) $ is $ \Gamma$--equivariantly homeomorphic to the Gromov boundary. This property is not true in general for Gromov hyperbolic spaces. Therefore we introduce

\begin{ass}
A Gromov hyperbolic space $ (X,d) $ satisifies the basic assumption (BA), if the horofunction boundary is $ \Isom(X,d)$--equivariantly homeomorphic to the Gromov boundary of $ ( X , d ) $.
\end{ass}
 
The assumption holds for all $ \CAT(-1) $--spaces, and by the remark above, for all hyperbolic groups with respect to a Green metric. 

We can now formulate a more general version of theorem \ref{thm:main}.

\begin{thm} \label{thm:main2}
Suppose $ (X,d) $ is a Gromov hyperbolic space which satisfies the geometric condition (BA). Suppose $ \mu $ is a symmetric probability measure on $ \Gamma $ with a finite
exponential moment, such that the group generated by the support of $ \mu $ acts non-elementary on $ \Xb $. Let $ g_1 , g_2, \ldots $ be a sequence of independent, $ \mu $--distributed random variables. Then, $ A > 0 $, and there is a positive constant $ \sigma $, such that
\[
\frac{1}{\sqrt{n}} \pan{d(g_1 \cdots g_n x_0,x_0) - nA} 
\]
converges weakly to a non-degenerate Gaussian distribution with variance $ \sigma $, as $ n \sra \infty $, and
\[
\limsup_{n \sra \infty} \frac{d(g_1 \cdots g_n x_0,x_0)-nA}{\sqrt{n \log{\log{n}}}} = \sigma > 0,
\]
almost surely with respect to $ \Pb $.
\end{thm}

For actions of linear semigroups on projective spaces, similar central limit theorems have been proved by Le Page \cite{Pa81}. A nice exposition can be found in \cite{BoLa85} and a more general
approach was later developed in \cite{HeHe04}. A standard and natural assumption in this theory is that the support of $ \mu $ is irreducible and ( strongly ) contracting ( \cite{BoLa85} ). This leads to a 
dynamical situation close to the action of a hyperbolic group, and lemma \ref{lem:spectralgap}, as its proof, closely follows the ideas in \cite{Pa81}.

This paper is organized as follows. Section \ref{section:rw} introduces the concepts of random walks and hyperbolic spaces. In subsection \ref{subsection:Green}, two different definitions
of the Green metric on a transient group are given, and we investigate some of its properties. In section \ref{section:ergodicrw}, we formulate the problem in the language of ergodic theory,
and derive a useful integral formula for the drift. In section \ref{section:final}, we present the complete proof of theorem \ref{thm:main2}.

\section{Random Walks on Hyperbolic Spaces} \label{section:rw}

\subsection{Random Walks on Isometry Groups}

Suppose $ (X,d) $ is a metric space. Let $ \Gamma $ be a subgroup of the isometry group $ \Isom(X,d) $, and let $ \mu $ be a probability measure on $ \Gamma $. Let
$ \mathfrak{B} $ denote the Borel $ \sigma $-algebra of $ \Gamma $, with respect to the compact--open topology. We assume that the measure $ \mu $ satisfies 
$ \int_{\Gamma} d(gx_0,x_0) \, d\mu(g) < +\infty $ for some point $ x_0 $ in $ X $. Define the probabilty measure space 
$ (\Omega,\mathfrak{F},\Pb) = (\Gamma^{\N}, \mathfrak{B}^{\N}, \mu^{\N})$. The $ \sigma $-algebra $ \mathfrak{F} $ has a natural filtration $ \mas{\mathfrak{F}_n}_{n \in \N}$ 
generated by cylinder sets in $ \Omega $.

We think of $ \Omega $ as the set of all measurable maps $ \N \sra \Omega $. The shift map $ T $ taking a function $ \omega \in \Omega $ to $ \omega(\cdot + 1) $ preserves the measure $ \Pb $ and is clearly ergodic. Let $ g $ be the projection $ \omega \mapsto \omega(0) $ from $\Omega $ onto $ \Gamma $. We will be interested in the asymptotic behaviour of the following ergodic cocycle into the group $ \Gamma $. Let $Z_0 = e $, where $ e $ is the identity in $ \Gamma $, and define for $ n \geq 1 $
\[
Z_{n}(\omega) = g(\omega) \cdots g(T^{n-1}\omega), \qquad \omega \in \Omega.
\]
Given some point $ x_0 $ in $ X $ we want to study the fluctuations of the random sequence $ A_n(\omega) = d(Z_n(\omega) x_0 , x_0 ) $, when $ n $ tends to infinity. 
It is obvious that $ A_n $ is an integrable subadditive cocycle, i.e. $ A_{n+m} \leq A_n + A_m \circ T^n $, for all $ n $ and $ m $, and hence the ergodic theorem of Kingman \cite{Ki73} 
asserts that the limit 
\[
A(\mu) = \lim_{n\sra \infty} \frac{A_n}{n} = \inf_{n \geq 1} \frac{1}{n} \int_{\Omega} A_n \, d \, \Pb =  \inf_{n \geq 1} \int_{\Gamma} d(\gamma x_0,x_0) \, d\mu^{*n}(\gamma) 
\]
exists and is constant a.e. with respect to $ \Pb $. We will refer to $ A(\mu) $ as the drift of the random walk $ Z_n $. Explicit calculations of $ A(\mu) $ are notoriously hard even in 
simple examples, see \cite{Gi07} for a nice exposition of formulae. In what follows, we will tactically assume we are in a situation when $ A(\mu) $ is positive.

\subsection{Horofunctions}

Let $ ( X, d ) $ be a metric space. A geodesic in $ X $ is an isometric embedding of an interval $ I \subset \R $ into $ X $.  We say that $ ( X, d ) $ is geodesic if any two points can
be connected by a geodesic. Note that we do not require the geodesic to be unique. A geodesic ray is a geodesic $ r $ defined on an interval $ [a,b) $ such that the limit of $ r(t) $
does not exist in $ X $ in the limit $ t \sra b^{-} $. 

A metric space $ ( X , d ) $ is proper if closed and bounded sets are compact. Pick a basepoint $ x_0 $ in $ X $ and define for every $ x \in X $ the function $ h_x(y) = d(x,y) - d(x,x_0) $. The 
map $ \Psi(x) = h_x $ is injective and the family $ F = \mas{h_x}_{x \in X} $ is equicontinuous in the space of continuous functions $ C(X) $. We define $ \X = \overline{F} $, where the closure of
$ F $ is taken in the topology of $ C(X) $. This is a compact Hausdorff space by the Arzela-Ascoli theorem, and it easy to check that any other choice of a basepoint $ x_0 $ gives rise to a homeomorphic
version of $ \X $. Furthermore, if $ ( X,d) $ is geodesic, the map $ x \sra h_x $ is a homeomorphism onto its image. Note that this is not necessarily true if the space is not geodesic. 

We denote the group of isometries on $ (X,d) $ by $ \Isom(X) $. We have a natural continuous action of $ \Isom(X) $ on $ \X $ by 
\[
(g.h)(x) = h(g^{-1}x) - h(g^{-1}x_0).
\]
We define the horofunction boundary of $ ( X , d ) $ to be the set $ \Xb = \X \backslash \Psi(X) $. It is clear that the action of $ \Isom(X,d) $ leaves $ \Xb $ invariant. Also note that $ \Xb $ is 
compact, if $ ( X , d ) $ is geodesic.

This compactification was promoted in Gromov's seminal paper \cite{Gr87} on hyperbolic groups, and its relevance in the theory of random walks has been demonstrated in series of papers
\cite{KaLe06}, \cite{KaLe07} and \cite{KaLe07_2}. 

\subsection{Gromov hyperbolic spaces}

Given three points $ x,y $ and $ z $ in $ X $, we define the Gromov product of $ x $ and $ y $ relative the point $ z $ by 
\[
(x,y)_z = \frac{1}{2}(d(x,z)+d(y,z)-d(x,y)).
\]
A metric space $ ( X,d) $ is Gromov hyperbolic if there is some $ \delta \geq 0 $, such that 
\[
(x,y)_w \geq \min\mas{(x,z)_w,(z,y)_w} - \delta,
\]
for all $ x,y,z $ and $ w $ in $ X $ ( It is an easy exercise to verify that $ \delta $ can be taken to be $ 0 $ for metric trees ). Examples of Gromov hyperbolic spaces include metric trees and manifolds 
of strictly negative sectional curvature. We say that a finitely generated group $ \Gamma $ is Gromov hyperbolic if its Cayley graph with respect to any finite generating set is a Gromov hyperbolic space.  
Examples include fundamental groups of compact negatively curved manifolds, e.g. surface groups of genus $ g \geq 2 $. It is a well-known fact that any Gromov hyperbolic group is finitely presented. In fact, most finitely presented groups are Gromov hyperbolic by a theorem by Olshanskii \cite{Ol92}:

\begin{thm}
Let $ A = \mas{a_1^{\pm 1},\ldots,a_k^{\pm 1}}$ be an alphabet, and let $ N(k,d,n) $ denote the number of group representations of the form $ \langle A \, | \, r_1 , \ldots r_d \rangle $, where $ \mas{r_i}_{i=1}^{d} $
is a set of reduced words of length $ n_i $ in the alphabet $ A $, and $ n = \pan{n_1,\ldots,n_d} \in \N^d $. Let $ N_h(k,d,n) $ denote the cardinality of the subset of hyperbolic groups in the above collection of groups. 
Then 
\[
N_h(k,d,n)/N(k,d,n) \sra 1 
\]
as $ \min\mas{n_1,\ldots,n_d} \sra \infty $, for fixed $ d $ and $ k $.
\end{thm}

Let $ ( X, d ) $ be a proper Gromov hyperbolic space. Let $ Y $ denote the space of all sequences which converge to infinity in the one-point compactification of $ X $. We say that two sequences in $ Y $
are equivalent if $ (x_i,y_j)_{x_0} \sra \infty $ as $ i,j \sra \infty $. This statement is clearly independent of the choice of the point $ x_0 $, and due to Gromov hyperbolicity, it is a transitive relation, and hence
an equivalence relation. Let $ \Xg $ be the set of all equivalence classes in $ Y $ under this relation. We can extend the Gromov product to $ \Xg $ by 
\[
(\xi,\eta)_{x_0} = \sup \liminf_{i,j \sra \infty} ( x_i,y_j)_{x_0},
\]
where $ x_i $ and $ y_i $ are sequences corresponding to $ \xi $ and $ \eta $, and the sup is taken over all such sequences. Note that $ (\xi,\xi)_{x_0} = \infty $. A natural candidate for a metric on $ \Xg $ 
would be $ \rho_{\epsilon}(\xi,\eta) = e^{-\eps (\xi,\eta)_{x_0}} $, for small $ \eps > 0 $, but unfortunately, this it not always a metric. However, it can always be deformed to a metric $ \rho $ on $ \Xg $, which satisfies, 
$ A\rho_\eps \leq \rho \leq B\rho_\eps $, for some constants $ A $ and $ B $ \cite{BrHa99}.

Suppose $ (X,d) $ is a Gromov-hyperbolic proper metric space. It is easy to see that there is a continuous $ \Isom(X,d) $--equivariant surjection from $ \Xb $ onto $ \Xg $ ( see \cite{St07} ), but 
this map is not always a homeomorphism. Indeed, let $ X $ be the Cayley graph of $ \Z \times \Z/{2\Z} $ with respect to the product of the standard generators on each group. The geometric realization is a bi-infinte ladder. Note that the horofunction boundary with respect to the word-metric ( extended to the edges ) is an interval, while the Gromov boundary only consists of two points. We will see below that this inconvenience disappears for transient hyperbolic groups, if we change 
from the word-metric to the quasi-isometric Green metric. 

\subsection{Green Metrics and Hilbert Metrics on Cones} \label{subsection:Green}

A measure $ \mu $ on $ \Gamma $ is symmetric if $ \mu(g^{-1)} = \mu(g) $ for all $ g \in \Gamma $, and finitely supported if $ \mu(g) \neq 0 $ for only a finite number of $ g $ in $ G $. 
Suppose $ \mu $ is a finitely supported and symmetric probability measure on a finitely generated group $ \Gamma $. Suppose that  the support of $ \mu $ generates $ \Gamma $ as a group. 
We define the Green kernel with respect to $ \mu $ as 
\[
G(x,y) = \sum_{n=0}^{\infty} \mu^{*n}(x^{-1}y)
\]
If the kernel $ G $ is everywhere finite for any finitely supported and symmetric measure, we say that $ \Gamma $ is transient. Recall the following result due to Varopoulos 
\cite{Va86}, which is a very strong 
generalization of Polya's classical theorem on random walks on $ \Z^d $.
\begin{thm}
A finitely generated group $ \Gamma $ is transient if and only if $ \Gamma $ is not virtually trivial, $ \Z $ or $ \Z^2 $.
\end{thm}

The Green metric $ d_G $ on a transient group $ \Gamma $ is defined as 
\[
d_G(x,y) = \log{\frac{G(e,e)}{G(x^{-1}y)}}, \qquad x,y \in \Gamma.
\]
We recall that the Green kernel on any transient group with respect to any finitely supported probability measure $ \mu $ can be factorized as $ G(x,y) = G(e,e)F(x,y) $, where $ F(x,y) $ denotes the 
probability for the $ \mu $--induced Markov chain to ever reach the point $ y $ in $ \Gamma $ from $ x $. Thus, the Green metric equals $ d_G(x,y) = - \log{F(x,y)} $. The metric axioms are easily verified under the above assumptions. The Green metric was introduced in \cite{BlBr07} in connection with DLA-processes on trees. Excellent accounts of applications of the Green metric to random walks and
the geometry of discrete groups can be found in \cite{Bl07} and \cite{Bl08}. This metric is also closely related to the Harnack metric, first introduced in \cite{Be65}.

\begin{exam}
Let $ \Gamma $ be the free group on $ q+1 $ generators. Let $ \mu $ be the uniform measure on the generators and the inverses. The Green function of $ \mu $ is given by 
\cite{Wo00}
\[
G(x,y) = \frac{q}{q-1} q^{-d'(x,y)},
\]
where $ d' $ denotes the word-metric with respect to the generators. Hence the Green metric $ d $ on $ \Gamma $ equals $ c d' $, where $ c = \log{q} $.
\end{exam}

The Green metric can also be introduced via Hilbert metrics. Let $ V $ be a topological vector space, equipped with an integrally closed partial order $ \preceq $ ( see \cite{Li95} ).
We define
\[
\alpha(f,g) = \sup \mas{\lambda > 0 \, | \, \lambda f \preceq g }
\]
\[
\beta(f,g) = \inf \mas{\lambda > 0 \, | \, g \preceq \lambda f  }
\]
We take $ \alpha = 0 $ and $ \beta = \infty $ if the corresponding sets are empty. We define the Hilbert ( pseudo-) metric on the projectivized cone $ P_\R \mathcal{C} $ as 
\[
d(f,g) = \frac{1}{2} \log{\frac{\beta(f,g)}{\alpha(f,g)}}
\]
The Hilbert metric associated to a cone in an integrally closed vector lattice has the universal property that any projective map between two cones is non-expansive with respect to the 
Hilbert metrics on the cones. 

For every $ y \in \Gamma $ we define the Martin kernel at $ y $ to be 
\[
K_y(x) = \frac{G(x,y)}{G(e,y)}.
\]  
It is easily verified that $ K_y $ is a superharmonic function on $ \Gamma $ with respect to the 
probability measure $ \mu $, i.e. $ K_y * \mu \leq K_y $ on $ \Gamma $. Furthermore, the map $ \Psi(y) = K_y $ is an embedding of $ \Gamma $ into the cone of positive superharmonic functions $ \mathcal{C} $ considered as a subset of the space of functions on $ \Gamma $ with the topology of pointwise convergence. It is proved in \cite{Wo00} that $ \mathcal{C} $ is a convex cone with a compact base. We define a $ \Gamma $-invariant metric $ \rho $ on $ \Gamma $ by restricting the Hilbert metric on this cone to the image of the embedding, i.e. 
\[
\rho(x,y) = d(K_x,K_y), \quad x,y \in \Gamma.
\]
We make the following simple observation which the author has not been able to locate in the literature.
\begin{lem}
Suppose $ \Gamma $ is a transient group. Then the induced metric $ \rho $ on the image $ \Psi(\Gamma) $ equals the Green metric on $ \Gamma $. 
\end{lem}

\begin{proof}
We first note that the inequality
\[
\frac{G(z,y)}{G(e,y)} \leq \lambda \frac{G(z,x)}{G(e,x)} \quad \forall \, z \in \Gamma,
\]
can be rewritten, in terms of the Green metric $ d $, as 
\[
e^{d(y,e)-d(x,e)} \sup_{z \in X} e^{d(z,x)-d(z,y)} \leq \lambda,
\]
and thus $ \lambda $ is finite. We conclude that 
\[
\beta(K_x,K_y) = e^{d(y,e)-d(x,e)} \sup_{z \in X} e^{d(z,x)-d(z,y)}  = e^{d(y,e)-d(x,e)} e^{d(y,x)} 
\]
The inequality for $ \alpha $ is completely analogous, 
\[
\alpha(K_x,K_y) = e^{d(y,e)-d(x,e)} \inf_{x \in X} e^{d(z,x)-d(z,y)} = e^{d(y,e)-d(x,e)} e^{-d(y,x)}.
\]
Hence, 
\[
\rho(K_x,K_y) = d(x,y),
\]
for all $x $ and $y $ in $\Gamma$. 
\end{proof}

\begin{rema}
It would be interesting to identify the projective endomorphisms of the cone $ \mathcal{C} $ which leave $ \Psi(\Gamma)$ invariant. These maps would induce natural 
semi--contractions on $ \Gamma $ with respect to the Green metric.
\end{rema}

\subsection{Connections between the Martin Boundary and the Gromov Boundary}

Suppose $ \Gamma $ is a discrete transient group, generated by a finite symmetric set $ S \subset \Gamma $. Let $ p $ be the uniform measure on $ S $. Then the Green
kernel with respect to $ p $ and the horofunction boundary with respect to the Green metric on $ \Gamma $ equals, by definition, the Martin boundary of $ \Gamma $ with 
respect to $ p $. In general, the Green metric is not equivalent to the word-metric with respect to $ S $. However, on transient 
hyperbolic groups, the two metrics are quasi-isometric ( see \cite{BlBr07} ), in the sense that there are constants $C$ and $b$, such that,
\[
\frac{1}{C}d_S(x,y) - b \leq d(x,y) \leq Cd_S(x,y) + b \qquad \forall \: x,y \in \Gamma,
\] 
where $d_S$ denotes the word--metric with respect to a finite generating set of $\Gamma$, and $d$ is the Green metric. In this case, we can let $b$ be zero in the inequality
on the right hand side. Note that the drift with respect to the Green metric of a 
random walk is positive if and only if it is positive with respect to the word--metric.

 We recall the following theorem by A. Ancona ( \cite{An87} and \cite{An90} ).

\begin{thm}
Suppose $ \Gamma $ is a transient Gromov hyperbolic group, and $ S $ is a finite generating set of $ \Gamma $. Let $ \mu $ be the uniform measure on $ S $, and let
$ d $ denote the Green metric with respect to the Green kernel of $ p $. Then the horofunction boundary of the metric space $ ( \Gamma, d ) $ is $ \Gamma $--equivariantly 
homeomorphic to the Gromov boundary of $ ( X , d ) $.
\end{thm}
In particular, for any transient Gromov hyperbolic group, $ ( \Gamma , d ) $ satisfies the geometric assumption (BA).

\subsection{General Aspects of Random Walks on Metric Spaces}

Let $ Y $ be a compact Hausdorff space with a continuous action of a group $ \Gamma $. If $ \mu $ is probability measure on $ \Gamma $ and $ \nu $ is a Borel probability  
measure on $ Y $, we define the convolution of $ \mu $ and $ \nu $ to be 
\[
(\mu * \nu )(\phib) = \int_{Y} \phib(gy) \, d\mu(g) d\nu(y), \quad \phib \in C(Y),
\]
which again is a Borel probability measure on $ Y $. If $ \nu $ is fixed, i.e. if $ \mu * \nu = \nu $, we say that $ \nu $ is $ \mu $-stationary, or simply stationary, if there is no
risk of confusion. Stationary measures always exist by a simple fixed point argument. The following theorem is due to V. Kaimanovich \cite{Ka00}

\begin{thm} \label{thm:Kaimanovich}
Suppose $ (X,d) $ is a Gromov hyperbolic space and $ \Gamma $ is a discontinuous subgroup of the isometry group of $ (X,d) $. Let $ \mu $ be a probability measure
on $ \Gamma $ such that the subgroup generated by the support of $ \mu $ is non-elementary with respect to $ \Xb $. Then the random walk on $ X $ induced by $ \mu $
converges almost everywhere in the hyperbolic compactification and there is a unique non-atomic $ \mu$-stationary measure $ \nu $ on $ \Xb $.
\end{thm}

Recall that an action of a hyperbolic group is non-elementary with respect to $ \Xb $ if the action does not fix any finite subset of $ \Xb $. 

We now turn to the question of positivity of the drift. For symmetric random walks on abelian and nilpotent groups it is easy to see that the drift is necessarily zero. However,
for non-amenable groups, Guivarc'h \cite{Gu80} 
proved the following theorem
\begin{thm} \label{thm:guivarch}
Let $ \Gamma $ be a non-amenable finitely generated group, and $ \mu $ a symmetric measure with finite first moment. Then $ A(\mu) > 0 $.
\end{thm}

It is a well-known fact that a hyperbolic group is amenable if and only if it is a finite group or virtually $ \Z $. 

In view of theorem \ref{thm:guivarch}, 
it is tempting to conjecture that Theorem \ref{thm:main}
is true for any symmetric random walk of finite support on a non-amenable group with 
respect to the word-metric. However, recall the following result by A. Erschler ( \cite{Dy99} and \cite{Er01} ), answering a question by A. Vershik, 
\begin{thm}
Let $ G = \Z \wr \Z $. There exists a symmetric probability measure $ \mu $ on $ G $ with finite support such that
\[
n^{-\frac{3}{4}}A_n(\mu), \sra 1 \quad n \sra \infty,
\]
where $ A_n(\mu) = \int_{G} d(g ,e) \, d \mu^{*n}(g) $, where $ d $ denotes some word-metric, and $ e $ is the identity element in $ G $.
\end{thm}
Let $ \Gamma = G \times \mathbb{F}_2 $, where $ \mathbb{F}_2 $ denotes the free group on two generators. Note that this group is not hyperbolic since the subgroup $\Z \wr \Z $
is an amenable non--hyperbolic group. Let $ \nu $ be the product of $ \mu $ and a finitely supported
symmetric measure on $ \mathbb{F}_2 $. The group $ \Gamma $ is non-amenable, but the random walk determined by the symmetric measure $ \nu $ will have linear drift 
with fluctuations of order $ n^{\frac{3}{4}} $. 

Counterexamples to the law of iterated logarithm can be constructed in a similar way using the results by D. Revelle \cite{Re03}. In particular, in \cite{Re03} it is proved that there 
are finitely supported measures on the group $ G $ above such that 
\[
0 < \limsup_{n \sra \infty} \frac{d(Z_n x_0,x_0) }{n^{3/4} (\log{\log{n}})^{1/4}} < \infty,
\]
almost surely. 

\section{Ergodic Theory of Random Walks on Proper Metric Spaces} \label{section:ergodicrw}

In this subsection we will derive a useful integral representation for the drift of a random walk generated by a symmetric probability measure $ \mu $. We assume that the
geometric assumption (BA) and the conditions in theorem \ref{thm:Kaimanovich} hold. Recall that $T$ denotes the forward shift map on $ \Omega = \Gamma^{\N} $, 
and let $ \tilde{\Omega} = \Omega \times H $. We note that the skew-product extension 
\[
\hat{T}(\omega,h) = ( T\omega , g(\omega)^{-1}.h ),
\]
has an invariant measure $\Pb = \Pb \times \nu$, where $\nu$ is the unique stationary measure on
the boundary $\partial \bar{X}$. 

Since $ \nu $ is non-atomic,
\[
A(\mu) = \lim_{n \sra \infty} \frac{d(Z_n x_0 , x_0)}{n} = \lim_{n \sra \infty} \frac{h(Z_n x_0)}{n},
\]
almost surely with respec to $\Pb$, for any fixed choice of $ h \in \Xg \cong \Xb $. Note that for any $h \in \Xb$, and $ \omega \in \Omega $, 
\[
h(Z_n(\omega) x_0) = Z_{n-1}^{-1}.h(g(T^{n-1}\omega)x_0) + h(Z_{n-1}(\omega)x_0) = \sum_{k=0}^{n-1}  Z_{k}^{-1}.h(g(T^k\omega)x_0),
\]
for $ n \geq 1 $. Thus, if we let $ F(\omega,h) = h(g(\omega)x_0) $, then 
\[
h(Z_n(\omega) x_0) = \sum_{k=0}^{n-1} F(\hat{T}^k(\omega,h)).
\]
It is obvious that $ F \in L^1(\Pbt) $, since $ \mu $ is assumed to satisfy $ \int_{\Gamma} d(gx_0,x_0) \, d\mu(g) $.  We have now outlined the main ingredients in the following important lemma.
\begin{lem} \label{lem:drift}
Suppose $ ( X , d ) $ is a Gromov hyperbolic space which satsifies the geometric assumption (BA), and $ \mu $ is a symmetric probability measure on $ \Gamma $ which 
satisfies the condtions in theorem \ref{thm:Kaimanovich}. Then, 
\[
A(\mu) = \int_{H} \int_{\Gamma} h(g x_0) \, d\mu(g) d\nu(h).
\]
\end{lem}
It should be noted that the above integral formula holds in a much more general context ( see \cite{KaLe08} ).

\section{Proofs of the theorems} \label{section:final}

The following lemma is a direct consequence of Gromov hyperbolicity and the basic assumption (BA).
\begin{lem} \label{lem:hypprop}
Suppose $ (X,d) $ is a Gromov hyperbolic space which satisfy the geometric property $ (BA) $. Suppose $ x_n \sra h' \in \Xb $, and $ h' \neq h $. Then the sequence 
\[
d(x_n,x_0) - h(x_n)
\]
is bounded.
\end{lem}

Suppose $ \mu $ is a symmetric probability measure on $ \Gamma $ which satisfies the conditions of Theorem \ref{thm:Kaimanovich}. Then, almost surely, the 
sequence $ Z_n x_0 $ converge to a point in the Gromov boundary $ \Xg \cong \Xb $, and for any Borel subset $ A \subset  \Xg $
\[
\Pb(\mas{\omega \in \Omega \: | \: Z_\infty(\omega) := \lim_{n \sra \infty} Z_n(\omega) x_0 \in A}) = \nu(A),
\]
where $ \nu $ is the unique stationary measure on $ \Xb $. Since $ \nu $ is non-atomic, and by lemma \ref{lem:hypprop},
\[
\frac{d(Z_n x_0 , x_0)-h(Z_n x_0)}{\sqrt{n}} \sra 0,
\]
almost surely $ \hak{\Pb} $, for any \emph{fixed} choice of $ h \in \Xb $. Thus, to prove theorem \ref{thm:main2},
to prove, for a fixed $ h \in \Xb $, 
\[
Y_n := \frac{h(Z_n x_0 , x_0 ) - nA(\mu)}{\sqrt{n}} \sra Y,
\]
in distribution, where $ Y $ is a centered and non-degenerate Gaussian variable.

Let us for a fixed $ h \in \Xb $ and $ u \in L^\infty(\Xb) $ define the sequence 
\[
M_n = h(Z_n x_0) - nA(\mu) + u(h) - u(Z_n^{-1}.h).
\]
We want to choose $ u $ so that $ M_n $ is a martingale with respect to the natural filtration $ \mas{\mathfrak{F}_n}_{n \geq 0} $. Since $ u $ is bounded, 
the sequence $ \frac{M_n}{\sqrt{n}} $ converge in distribution to a non-degenerate and centered Gaussian variable $ Y $ if and only if the sequence $ Y_n $ 
converge in distribution to $ Y $.

Note that 
\[
\mathbb{E}\hak{h(Z_n(\cdot) g(T^{n} \cdot )x_0,x_0) - h(Z_n(\cdot)x_0,x_0) \: | \: \mathfrak{F}_n } =  \int_{\Gamma} Z_{n}(\cdot)^{-1}.h(g x_0) \, d\mu(g), 
\] 
almost everywhere $ \hak{\Pb} $. Thus, if we introduce the operator 
\[
P \phib(h) = \int_{\Xb} \phib(g^{-1}.h) \, d\mu(g), \quad \phib \in L^{\infty}(\Xb,\nu), 
\]
we see that if we can choose $ u \in L^\infty(\Xb,\nu) $ such that
\[
( I - P)u = \psi,
\]
where $ \psi(h) = \int_{\Gamma} h(gx_0) \, d\mu(g) - A(\mu) $, then $ M_n $ is a martingale sequence with respect to the filtration $  \mas{\mathfrak{F}_n}_{n \geq 0}  $.

It is clear that the equation above makes sense in any normed space of bounded functions on $ \Xb $. However, to guarantee the existence of a solution we need to 
restrict the operator to a sufficiently nice space of functions on $ \Xb $. We observe that 
\[
\int_{\Xb} \psi(h) \, d\nu(h) = \int_{\Xb} \int_{\Gamma} h(gx_0) \, d\mu(g) d\nu(h) - A(\mu) = 0, 
\]
by lemma \ref{lem:drift}. 

Let $ L^1_0(\Xb) $ denote the space of all $ \nu $--integrable functions with integral zero. Since $ \mu $ is symmetric, $ P $ preserves $ L^1_0(\Xb) $. Recall that a real valued 
function $ \phib $ on a metric space $ ( Y , \rho ) $ is H\"older continuous if 
\[
\sup_{y \neq z} \frac{|\phib(y)-\phib(z)|}{\rho(y,z)^{\alpha}} < \infty,
\]
for some $ \alpha > 0 $. The space of H\"older continuous functions $ \H_\alpha $ on $ ( \Xb , \rho ) $ for a fixed $ \alpha > 0 $ is a Banach space with respect to the norm
\[
||\phib||_\alpha = | \int_{\Xb} \phib(h) \, d\nu(h) | + \sup_{h \neq h'} \frac{|\phib(h)-\phib(h')|}{\rho(h,h')^{\alpha}}.
\]
Note that $ \H_\alpha^0 = \H_\alpha \cap L^1_0(\Xb,\nu) $ is a closed subspace of $ \H_\alpha $. We want to solve the equation $ ( I - P )u = \psi $ for $ u \in \H_\alpha^0 $,
for some $ \alpha > 0 $. To guarantee the existence of a solution to this equation in the above space, it suffices to prove that the spectral radius of $ P $ in $ \H_\alpha^0 $ is
less than $ 1 $, at least for sufficiently small $ \alpha > 0 $.

We say that a probability measure on $ \Gamma $ has exponential moment if, for some $ \beta > 0 $, the integral $ \int_{\Gamma} e^{\beta d(g x_0 , x_0 ) } \, d\mu(g) $  is finite.
The following lemma is the key result of this section.
\begin{lem} \label{lem:spectralgap}
Suppose $ \mu $ has a finite exponential moment, and that $ ( X , d) $ satisfies the geometric assumption (BA). Then $ ||P^n||_\alpha \leq C \tau^n $, for some positive constant $ C $, 
for all $ n \geq 1 $ and $ \alpha $ sufficiently small. 
\end{lem}

\begin{proof}
It is easy to see that it suffices the following strong proximality statement on the boundary $ \partial X $,
\[
\sup_{h \neq h' } \int_{\Gamma} \Big(\frac{\rho(g^{-1} h,g^{-1} h')}{\rho(h,h')} \Big)^{\alpha} \, d\mu^{*n}(g) < 1,
\]
for sufficiently small $ \alpha > 0 $ and $ n $ large. This estimate is equivalent to ( we identify a point in $ \Xg $ with its horofunction ),
\[
\sup_{h \neq h'} \int_{\Gamma} e^{\alpha((h,h')_{x_0} - (g^{-1}h,g^{-1}h')_{x_0})} \, d \mu^{*n}(g) < 1,
\]
since $ \rho(\cdot,\cdot) $ is uniformly equivalent to $ e^{-(\cdot,\cdot)_{x_0}} $. For small $ \alpha $, and large $ n $, we can, due to the finite 
exponential moment condition on $\mu$, replace the integrand ( up to an arbitrary small error which depends on $\alpha$ ) with
\[ 
1 - \alpha((g^{-1}h,g^{-1}h')_{x_0} - (h,h')_{x_0}).
\]
Because of assumption (BA), we have the identity, 
\[
 (h,h')_{x_0} - (g^{-1}h,g^{-1}h')_{x_0} = 2( h(gx_0) + h'(gx_0) ),
\]
where $\rho$ denotes the metric on the boundary. Thus, it suffices to prove that
\[
\inf_{h \in \Xb} \frac{1}{n} \int_{\Gamma} h(gx_0) \, d\mu^{*n}(g) = \inf_{h \in \Xb} \frac{1}{n} \int_{\Omega} h(Z_n x_0) \, d\Pb > 0.
\]
Note that the order of $ \inf $ and the integral is crucial here. Suppose that the last equality does not hold. 
Then we can find a sequence of integers $ n_k $, tending to infinity, and horofunctions $ h_k $ such that 
\[ 
\liminf_{k \sra \infty} \frac{1}{n_k} \int_{\Omega} h_k(Z_{n_k}x_0) \, d\Pb \leq 0,
\]
and $ h_k $ converges in $ \X $ to some $ h $. This implies that the sequence  $ Z_{n_k}x_0 \sra h $ in $ \Xb $, with a positive probability, which is a 
contradiction to theorem \ref{thm:Kaimanovich}, since the stationary distribution on $ \Xb $ is non-atomic. 
\end{proof}

\begin{rema}
In order to deal with the case of non--symmetric measures, we introduce the operator $N$ which performs integration with
respect to the unique stationary measure $\check{\nu}$ on $\partial \overline{X}(\infty)$ with respect to the measure $\check{\mu}(g) = \mu(g^{-1})$, 
and extend the definition of $P$ to $\mathcal{H}_\alpha$. A slight modification of the proof of Lemma \ref{lem:spectralgap} will give that 
$||P^n-N||_\alpha$ decays to zero exponentially fast as $n \sra \infty$. The arguments above ( for a general measure $\mu$, not necessarily symmetric )
now lead to the equation $(I-P)u = (I-N)\psi$, which can be solved for $u$ in $\mathcal{H}_\alpha$ since $ \sum_{n \geq 0} P^n (I-N)\psi $ converges in 
$\mathcal{H}_\alpha$ by the exponential decay of $||P^n - N||_{\alpha}$ and the fact that $P^n N = N$ for all $n \geq 1$. A more detailed 
description in a similar situation can be found in the paper \cite{Le01}.
\end{rema}

Lemma \ref{lem:spectralgap} is related to similar spectral gap results achieved by Le Page \cite{Pa81} and Guivarc'h \cite{Gu80}, \cite{GuPa04} for actions of linear subgroups on projective spaces. A standard assumption in this setting is that the semigroup generated by the support of $ \mu $ is contracting and irreducible. This leads to a dynamical situation close to actions of hyperbolic groups on the Gromov boundary.

This approach closely follows the ideas on martingale approximations developed by M. Gordin in \cite{Go69}. Very nice and general treatments of martingale approximations, 
quasi-compact operators and connections to central limit theorems can be found in \cite{HeHe04} and \cite{GoHo04}. 

The following Lindeberg--type central limit theorem and law of iterated logarithm for martingales can be found in \cite{HaHe80} ( section 1.7 ).
\begin{thm}
Suppose $ (\Omega,\mathfrak,\Pb) $ be a probability space, and assume that $ \mas{\mathfrak{F}_n}_{n} $ is a filtration of $ \mathfrak{F} $.
Let $ M_n $ be a centered martingale, with respect to this filtration, and set $ X_k = M_k - M_{k-1} $, for $ k \geq 1 $, and $ M_0 = 0 $. 
Suppose that, $ M_n $ is a $ L^p $--martingale for all $ p < \infty $, and for all $ \eps > 0 $, 
\[
\frac{1}{n} \sum_{k = 1}^{n} \mathbb{E}\hak{X_{k}^2 \, \chi_{|X_k| > \eps \sqrt{n}} \, | \, \mathfrak{F}_{k-1} } \sra 0,
\]
almost surely, and
\[
\frac{1}{n} \sum_{k=1}^{n} \mathbb{E}\hak{X_k^2 \, | \, \mathfrak{F}_{k-1} } \sra \sigma^2,
\]
almost surely, where $ \sigma $ is an essentially bounded measurable function, and $ \mathbb{E} $ denote the conditional expectation operator with respect
to $ ( \mas{\mathfrak{F}_{k}}_{k \geq 1} , \Pb ) $. Then
\[
\frac{M_n}{\sqrt{n}} \sra M,
\]
in distribution, where the distribution $ \mu $ of $ M $ has the Fourier transform $\hat{\mu}(\omega)=\mathbb{E}\hak{\exp{-\frac{1}{2}\eta^2\omega^2}}$. In 
particular, if $ \sigma $ is essentially constant, $ M $ is Gaussian. In this case, we also have
\[
-\sigma = \liminf_{n \sra \infty} \frac{M_n}{\sqrt{n \log{\log{n}}}} \quad \textrm{and} \quad \limsup_{n \sra \infty} \frac{M_n}{\sqrt{n \log{\log{n}}}} = \sigma,
\]
almost surely.
\end{thm}

In our case, 
\begin{eqnarray*}
\mathbb{E}\hak{X_k^2 \, | \, \mathfrak{F}_k} & = & \mathbb{E}\hak{ \big(h(Z_k x_0) - h(Z_{k-1}x_0) - A(\mu) + u(Z_{k}^{-1}.h) - u(Z_{k-1}^{-1}.h) \big)^2 \, | \, \mathfrak{F}_k}  \\
& = & \int_{\Gamma} \big( Z_{k-1}^{-1}.h(g x_0) - A(\mu) + u(g^{-1} Z_{k-1}^{-1}.h) - u(Z_{k-1}^{-1}.h) \big)^2  \, d\mu(g) \\
& = & G(\hat{T}^{k-1}(\cdot,h)),
\end{eqnarray*}
where
\[
G(\omega,h) = \int_{\Gamma} \big( h(gx_0) - A(\mu) + u(g^{-1}.h) - u(h) \big)^2 \, d\mu(g).
\]
in the notation of subsection \ref{section:ergodicrw}. Since $ \mu $ is assumed to have a finite exponential moment, $ G \in L^p(\Pbt) $, for all $ 1 \leq p < \infty $, and thus
the ergodic theorem apply
\begin{eqnarray*}
\sigma^2 & = & \lim_{n \sra \infty} \frac{1}{n} \sum_{k=1}^{n} \mathbb{E}\hak{X_k^2 \, | \, \mathfrak{F}_k} = \int_{\Omega} \int_{H} G(\omega,h) \, d\Pb(\omega) d\nu(h) \\
& = & \int_{H} \int_{\Gamma} \big( h(gx_0) - A(\mu) + u(g^{-1}.h) - u(h) \big)^2 \, d\mu(g) d\nu(h).
\end{eqnarray*}
In particular, $ \sigma $ is essentially constant, and $ \frac{M_n}{\sqrt{n} } $ converge in distribution to a centered Gaussian variable, which a priori can be degenerate, i.e. $ \sigma = 0 $.
However, this would entail that
\[
h(gx_0) = A(\mu) + u(g^{-1}.h) - u(h) 
\]
for all $ g \in \mathrm{supp}{(\mu)} $ and for almost every $ h $ in $ \Xb $. Thus, if $ g_1 $ and $ g_2 $ are in $ \mathrm{supp}(\mu) $, and $ h $ is a $ \nu $--generic point in $ \Xb $, then
\[
h(g_1g_2x_0) = g_1^{-1}.h(g_2 x_0) + h(g_1 x_0) = 2A(\mu) + u((g_1 g_2)^{-1}.h) - u(h),
\]
and so we can conclude that for any non-trivial element $ g $ in the support of $ \mu $, there is some positive integer $ c(g) $, such that
\[
h(gx_0) = c(g) A(\mu) + u(g^{-1}.h) - u(h) 
\]
for every $ \nu $--generic point in $ \Xb $. However, since the support was assumed to generate a non-elementary group $ \Gamma $, there must be a hyperbolic element $ g $ in the support
of $ \mu $, such that its fixed points $ h $ and $ h' $ in $ \Xb $ are generic for the measure $ \nu $. Thus, by proximality of the action on the boundary, 
\begin{eqnarray*}
h(gx_0) & = & \lim_{n \sra \infty} d(gx_0,g^{n}x_0) - d(g^{n}x_0,x_0) \\
& = & - \lim_{n \sra \infty} d(gx_0,g^{-n}x_0) - d(g^{-n}x_0,x_0) = -h'(gx_0),
\end{eqnarray*}
and we get the equation 
\[
h(gx_0) = A = -h(gx_0),
\]
which implies $ A = 0 $, and therefore contradicts Theorem \ref{thm:guivarch}. We have thus proved the following theorem.
\begin{thm}
Suppose $ (X,d) $ is a Gromov hyperbolic space which satisfies the geometric condition (BA). Suppose $ \mu $ is a symmetric probability measure on $ \Gamma $ with an exponential moment,
such that the group generated by the support of $ \mu $ acts non-elementary on $ \Xg $. Then there is a positive constant $ \sigma $, such that
\[
\frac{1}{\sqrt{n}} \pan{d(Z_n x_0,x_0) - nA(\mu)} 
\]
converge weakly to a non-degenerate Gaussian distribution, as $ n \sra \infty $, and
\[
\limsup_{n \sra \infty} \frac{d(Z_n x_0,x_0)-nA(\mu)}{\sqrt{n \log{\log{n}}}} = \sigma > 0,
\]
almost surely with respect to $ \Pb $.
\end{thm}

The necessary extension to cover the case of non-symmetric measures is standard. See \cite{Le01}.

\section{Acknowledgements}
The author is grateful to Anders Karlsson and Uri Bader for many interesting discussions on this paper.

\end{document}